\documentclass{amsart}

\usepackage{amssymb, amsmath}

\newcommand{\N}{\mathbb{N}}
\newcommand{\Z}{\mathbb{Z}}

\newcommand{\Q}{\mathbb{Q}}
\newcommand{\R}{\mathbb{R}}
\newcommand{\CC}{\mathbb{C}}

\newcommand{\goth}{\mathfrak}
\newcommand{\mc}{\mathcal}
\newcommand{\ord}{\mathop{\mathrm{ord}}\nolimits}

\newtheorem{theorem}{Theorem}[section]
\newtheorem{lemma}[theorem]{Lemma}
\newtheorem{corollary}[theorem]{Corollary}

\theoremstyle{remark}

\numberwithin{equation}{section}
\begin{document}

\title{New Spectral Multiplicities for Ergodic Actions}

\author[A. V. Solomko]{Anton V. Solomko}
\address{Institute for Low Temperature Physics \& Engineering of National Academy of Sciences of Ukraine, 47 Lenin Ave., Kharkov, 61164, UKRAINE}
\email{solomko.anton@gmail.com}


\begin{abstract}
Let $G$ be a locally compact second countable Abelian group.
Given a measure preserving action $T$ of $G$ on a standard probability space $(X, \mu)$,
let $\mc M(T)$ denote the set of essential values of the spectral multiplicity function of the Koopman representation $U_T$ of $G$ defined in $L^2(X,\mu)\ominus \CC$ by $U_T(g)f := f\circ T_{-g}$.
In the case when $G$ is either a discrete countable Abelian group or $\R^n$, $n\geqslant 1$, it is shown that the sets of the form $\{p,q,pq\}$, $\{p,q,r,pq,pr,qr,pqr\}$ etc. or any multiplicative (and additive) subsemigroup of $\N$ are realizable as $\mc M(T)$ for a weakly mixing $G$-action $T$.
\end{abstract}

\maketitle

\section{Introduction}

Let $G$ be a locally compact second countable Abelian group and let $T=(T_g)_{g\in G}$ be a measure preserving action of $G$ on a standard probability space $(X, \goth B, \mu)$.
Denote by $U_T$ the induced Koopman unitary representation of $G$ in $L^2_0(X,\mu) := L^2(X,\mu)\ominus \CC$ given by
$$
U_T(g)f := f\circ T_{-g}.
$$
By the spectral theorem, there is a probability measure $\sigma$ on the dual group $\widehat{G}$ called a \emph{measure of maximal spectral type} of $U_T$ and a measurable field of Hilbert spaces $\widehat{G}\ni\omega \mapsto \mc H_\omega$ such that
$$
L^2_0(X,\mu) = \int_{\widehat{G}}^\oplus \mc H_\omega d\sigma(\omega)
\quad \text{and}\quad
U_T(g) = \int_{\widehat{G}}^\oplus \omega(g)I_\omega d\sigma(\omega), \; g\in G,
$$
where $I_\omega$ is the identity operator on $\mc H_\omega$ \cite{Nai}.
A map $m_T\colon\widehat{G}\ni\omega \mapsto \dim\mc H_\omega \in \N\cup\{\infty\}$ is called the \emph{spectral multiplicity function} of $U_T$.
Let $\mc M(T)$ stand for the set of essential values of $m_T$.
We are interested in the following \emph{spectral multiplicity problem}:
\begin{enumerate}
\item[({\bf Pr})]
Which subsets $E \subset \N$ are realizable as $E = \mc M(T)$ for an ergodic (or weakly mixing) $G$-action $T$?
\end{enumerate}

This problem was studied by a number of authors (see the recent survey \cite{Da_surv} and references therein) mainly in the case $G=\Z$.
It is proved, in particular, that a subset $E\subset\N$ is realizable in each of the following cases:
\begin{itemize}
\item
$1\in E$ (\cite{KwL} for $G=\Z$, \cite{DL} for $G=\R$),
\item
$2\in E$ (\cite{KaL} for $G=\Z$, \cite{DL} for $G=\R$),
\item
$E=\{p\}$ for arbitrary $p\in\N$ (\cite{Ag}, \cite{Ry}, \cite{Da_hom} for $G=\Z$, \cite{DS} for $\R^n$ and arbitrary discrete countable Abelian group),
\item
$E = n\cdot F$ for arbitrary $F\ni 1$ and $n>1$ (\cite{Da_hom} for $G=\Z$).
\end{itemize}

Our aim is to obtain some new spectral multiplicities first appeared in \cite{Ry} for $G=\Z$.
Given $E,F\subset \N$, let $E\diamond F := E \cup F \cup EF$\footnote{Given $E,F\subset \N$, by $EF$ we mean their algebraic product, i.e. $EF=\{ef \mid e\in E, f\in F\}$.}.
In this notation, $\{p\} \diamond \{q\} = \{p,q,pq\}$, $\{p\} \diamond \{q\} \diamond \{r\} = \{p,q,r,pq,pr,qr,pqr\}$ etc.

\begin{theorem}\label{MainTh}
Let $G$ be either a discrete countable Abelian group or $\R^m$ with $m\geqslant 1$.
Given a (finite or infinite) sequence of positive integers $p_1,p_2,\ldots$, there exists a weakly mixing probability measure preserving $G$-action $T$ such that $\mathcal M(T) = \{p_1\}\diamond\{p_2\}\diamond\cdots$.
\end{theorem}

Since any multiplicative subsemigroup of $\N$ can be represented in the form $\{p_1\}\diamond\{p_2\}\diamond\cdots$, we obtain the following

\begin{corollary}
Any multiplicative (and hence any additive) subsemigroup $E$ of $\N$ is realizable as $E=\mc M(T)$ for a weakly mixing $G$-action $T$.
\end{corollary}

To prove Theorem~\ref{MainTh} we adapt the idea from \cite{Ry}.
The required action is the product $T_1\times T_2\times\cdots$, where $T_i$ is a weakly mixing $G$-action with homogeneous spectrum of multiplicity $p_i$.
The existence of such actions was proved in \cite{DS} via `generic' argument originated from \cite{Ag}.
To `control' the spectral multiplicities of Cartesian products of such actions we furnish $T_i$
with certain asymptotical operator properties 
using both `generic' argument and $(C,F)$-technique.

In Section~\ref{Prelim_sec} we list some basic definitions and facts that will be used in the sequel to prove the main theorem.
Only Subsection~\ref{Unitary_sec} contains the detailed proofs of some original results related to the spectral multiplicities for unitary representations.
In Subsection~\ref{CF_sec} we briefly outline the $(C,F)$-construction of measure preserving actions which is an algebraic counterpart of the classical geometric `cutting-and-stacking' technique and in \ref{PoisSusp_sec} we recall the definition and some basic properties of the Poisson suspension that allows us to obtain finite measure preserving actions from infinite measure preserving ones.
Both techniques are used to construct explicitly rigid actions in Lemmata~\ref{LemRigid} and \ref{LemDIlim}.
In Section~\ref{Rp_sec} we prove Theorem~\ref{MainTh} in the case where $G = \R^m$.
In general, the proof goes along the lines developed in \cite{Ry}.
To prove Theorem~\ref{MainTh} for arbitrary discrete countable Abelian group we need some modification of this scheme.
This is done in Section~\ref{Disc_sec}.
Though both profs can be given in spirit of Section~\ref{Disc_sec}, the constraints appeared in Section~\ref{Disc_sec} seem to be artificial and this is the main reason why we consider separately two cases for $G$.

\subsection*{Acknowledgement}
I would like to thank my advisor Alexandre Danilenko for introducing this subject to me and for numerous helpful discussions.

\section{Preliminaries} \label{Prelim_sec}



\subsection{Unitary representations} \label{Unitary_sec}

Denote by $\mc U(\mc H)$ the group of unitary operators on a separable Hilbert space $\mc H$.
We endow $\mc U(\mc H)$ with the (Polish) strong operator topology (which on $\mc U(\mc H)$ is also the weak operator topology).
Given a locally compact second countable group $\Gamma$, we furnish the product space $\mc U(\mc H)^\Gamma$ with the (Polish) topology of uniform convergence on the compact subsets in $\Gamma$.
Denote by $\mc U_\Gamma(\mc H)\subset\mc U(\mc H)^\Gamma$ the subset of all unitary representations of $\Gamma$ in $\mc H$.
Obviously, $\mc U_\Gamma(\mc H)$ is closed in $\mc U(\mc H)^\Gamma$ and hence Polish in the induced topology.
Let $\mc B(\mc H)$ stand for the set of all boundary linear operators on $\mc H$ endowed with the weak operator topology.
By a \emph{unitary polynomial on $\Gamma$} we mean a mapping $P\colon \mc U_\Gamma(\mc H) \to \mc B(\mc H)$ in the form
$$
P(U) = \alpha_1 U(g_1) + \cdots + \alpha_n U(g_n), \quad
\alpha_i \in \CC, g_i \in \Gamma, U \in \mc U_\Gamma(\mc H).
$$
We now list some lemmata that will be needed while proving the main theorem.

\begin{lemma} \label{Pol}
Given a unitary polynomial $P\colon \mc U_\Gamma(\mc H) \to \mc B(\mc H)$ and a sequence $(g_n)_{n=1}^\infty$ in $\Gamma$, the set
$$\mc P := \{ U \in U_\Gamma(\mc H) \mid P(U) \text{ is a limit point of } \{U(g_n)\}_{n\in \N} \}$$
is $G_\delta$-subset in $U_\Gamma(\mc H)$.
\end{lemma}

\begin{proof}
Let $d$ stand for a metric compatible with the week topology on $\mc B(\mc H)$.
Then
$$
\mc P = \bigcup_{m=1}^\infty \bigcup_{N=1}^\infty \bigcap_{n=N}^\infty \{ U \in U_\Gamma(\mc H) \mid d(P(U),U(g_n)) < \frac{1}{m} \}.
$$
Obviously, the sets $\{ U \in U_\Gamma(\mc H) \mid d(P(U),U(g_n)) < \frac{1}{m} \}$ are open in $U_\Gamma(\mc H)$.
\end{proof}

Recall that two unitary representations $U,V\in \mc U_G(\mc H)$ of an Abelian group $G$ are called \emph{spectrally disjoint} if their measures of maximal spectral type $\sigma_U$ and $\sigma_V$ are mutually singular: $\sigma_U \perp \sigma_V$.
By $\mc M(U)$ we denote the essential image of the spectral multiplicity function of $U$.
It is clear that if $U$ and $V$ are spectrally disjoint then $\mc M(U\oplus V) = \mc M(U) \cup \mc M(V)$.
Lemma~\ref{LemDis} gives us the useful sufficient condition of spectrally disjointness.

\begin{lemma} \label{LemDis}
Let $G$ be a locally compact second countable Abelian group.
Let $U,V\in \mc U_G(\mc H)$.
If there is a sequence $(g_n)_{n=1}^\infty \subset G$ such that
$$
U(g_n)\to I \text{ and } V(g_n)\to 0,
$$
then $U$ and $V$ are spectrally disjoint.
\end{lemma}

\begin{proof}
Let $\sigma_U$ and $\sigma_V$ be measures of maximal spectral type of $U$ and $V$ respectively.
By the spectral theorem,
$$
\mc H = \int_{\widehat{G}}^\oplus \mc H_\omega^{(1)} d\sigma_U(\omega), \quad
U(g)=\int_{\widehat{G}}^\oplus \omega(g)I_\omega d\sigma_U(\omega),
$$
$$
\mc H = \int_{\widehat{G}}^\oplus \mc H_\omega^{(2)} d\sigma_V(\omega), \quad
V(g)=\int_{\widehat{G}}^\oplus \omega(g)I_\omega d\sigma_V(\omega).
$$
Suppose $\sigma_U$ is equivalent to $\sigma_V$ on some subset $A\subset \widehat{G}$ with $\sigma_U(A)>0$.
Take any $0\neq f \in \mc H$ with $\text{supp} f \subset A$.
Then on the one hand
$$
\int_{\widehat{G}} \langle\omega(g_n)f(\omega),f(\omega)\rangle d\sigma_U(\omega) =
\langle U(g_n)f,f \rangle \to \|f\|^2 \neq 0.
$$
On the other hand
\begin{multline*}
\int_{\widehat{G}} \langle\omega(g_n)f(\omega),f(\omega)\rangle d\sigma_U(\omega) =\\
= \int_{\widehat{G}} \langle\omega(g_n)f(\omega),f(\omega)\rangle \frac{d\sigma_U}{d\sigma_V}(\omega) d\sigma_V(\omega) =
\langle V(g_n)f, \frac{d\sigma_U}{d\sigma_V} f \rangle \to 0.
\end{multline*}
This contradiction proves that $\sigma_U \perp \sigma_V$.
\end{proof}

Given  $U,V \in \mc U_G(\mc H)$, by their \emph{tensor product} we mean the unitary representation $U \otimes V$ of $G$ in $\mc H \otimes \mc H$ defined by $(U \otimes V) (g) := U(g) \otimes V(g)$.
If $\sigma_U$ and $\sigma_V$ are measures of maximal spectral type of $U$ and $V$, then the convolution $\sigma_U \ast \sigma_V$ is a measure of maximal spectral type of $U \otimes V$.
Let
$$
\sigma_U \times \sigma_V = \int_{\widehat{G}} \sigma_\omega d (\sigma_U \ast \sigma_V) (\omega)
$$
stand for the disintegration of $\sigma_U \times \sigma_V$ with respect to the projection map $\widehat{G}\times\widehat{G} \ni (\omega_1,\omega_2) \mapsto \omega_1 \omega_2 \in \widehat{G}$.
Then the map $\widehat{G} \ni \omega \mapsto \dim( L^2(\widehat{G}\times\widehat{G}, \sigma_\omega) )$ is the multiplicity function of $U \otimes V$.

The following lemma which is an obvious generalization of \cite[Lemma~3.1]{Ry} allows us to `control' the spectral multiplicities of tensor products.
Recall that a unitary representation $U \in \mc U_G(\mc H)$ has \emph{simple spectrum} (i.e. $\mc M(U) = \{1\}$) if and only if there is $\varphi \in \mc H$ (called a \emph{cyclic vector} for $U$) such that the smallest closed subspace $\mc H_\varphi$ of $\mc H$ containing all the vectors $U(g)\varphi$, $g\in G$, is the entire $\mc H$.
$\mc H_\varphi$ is called  the \emph{cyclic subspace} of $\varphi$.

\begin{lemma} \label{LemRSS}
Let $G$ be a locally compact second countable Abelian group and let $U,V\in \mc U_G(\mc H)$.
Suppose there exists a sequence $(g_n)_{n=1}^\infty \subset G$ and its subsequences $(g_{n_k(i)})_{k=1}^\infty$, $i \in J$, such that
\begin{enumerate}
\item \label{limi}
$U(g_n)\to I$ as $n\to\infty$ and
\item \label{limii}
$V(g_{n_k (i)})\to V(d_i)$ as $k\to\infty$ for each $i \in J$,
\end{enumerate}
where $\{ d_i \}_{i \in J} \subset G$ is at most countable subset such that $\langle d_i \rangle_{i \in J}$\footnote{Given a subset $A \subset G$, by $\langle A \rangle$ we denote the smallest subgroup of $G$ containing $A$.} is dense in $G$.
Then
\begin{enumerate}
\item[(1)]
if $U$ and $V$ have simple spectrum then $U\otimes V$ has simple spectrum;
\item[(2)]
in general case, $\mc M(U\otimes V) = \mc M(U) \mc M(V)$.
\end{enumerate}
\end{lemma}

\begin{proof}
(1)
Let $\varphi$ and $\psi$ be cyclic vectors for $U$ and $V$ respectively.
We claim that $\varphi\otimes\psi$ is a cyclic vector for $U\otimes V$.
Indeed, the cyclic subspace $\mc H_{\varphi\otimes\psi}$ of $\varphi\otimes\psi$ is weakly closed\footnote{Here we use the fact that any (strongly) closed convex set is weakly closed.}, invariant under $U(g)\otimes V(g)$ for each $g\in G$ and contains all the vectors $U(g)\varphi \otimes V(g)\psi$, $g\in G$.
Hence by (\ref{limi}) and (\ref{limii}) it contains all the weak limits
\begin{align*}
\varphi\otimes V(d_i)\psi & = \lim_{k \to\infty} U(g_k(i))\varphi \otimes V(g_k(i))\psi, \\
\varphi\otimes V(d_i + d_j)\psi & = \lim_{k \to\infty} U(g_k(j))\varphi \otimes V(g_k(j))V(d_i)\psi, \\
\text{etc.} &
\end{align*}
The space $\mc H_{\varphi\otimes\psi}$ contains therefore all the vectors $\varphi\otimes V(d)\psi$, $d\in \langle d_i \rangle_{i\in J}$.
Since $\mc H_{\varphi\otimes\psi}$ is invariant under $U(g)\otimes V(g)$ for each $g\in G$ it contains all the vectors $U(g)\varphi\otimes V(d+g)\psi$, $g\in G$, $d\in \langle d_i \rangle_{i\in J}$, which form a total system in $\mc H\otimes\mc H$.
Hence $U\otimes V$ has simple spectrum.

(2)
Let
$$
U = \bigoplus_{p\in \mc M(U)} p U^{(p)} \quad \text{and} \quad V = \bigoplus_{q\in \mc M(V)} q V^{(q)},
$$
where $U^{(p)}$ (and $V^{(q)}$) are spectrally disjoint and have simple spectrum.
In other words, $\bigoplus_p U^{(p)}$ and $\bigoplus_q V^{(q)}$ have simple spectrum.
Then for $U \otimes V$ we have the following decomposition:
$$
U \otimes V = \bigoplus_{\substack{p\in \mc M(U)\\q\in \mc M(V)}} pq (U^{(p)} \otimes V^{(q)}).
$$
As we have already shown in (1), $\bigoplus_{p,q} U^{(p)} \otimes V^{(q)} = \bigoplus_p U^{(p)} \otimes \bigoplus_q V^{(q)}$ has simple spectrum.
This means that $U^{(p)} \otimes V^{(q)}$, $(p,q)\in \mc M(U)\times\mc M(V)$, are spectrally disjoint and have simple spectrum.
Hence $\mc M(U\otimes V) = \mc M(U) \mc M(V)$.
\end{proof}

Following \cite{Ry}, we will say that $U$ and $V$ are \emph{strongly disjoint} if the map $(\widehat{G}\times\widehat{G},\sigma_U \times\sigma_V) \ni (\omega_1,\omega_2) \mapsto \omega_1 \omega_2 \in (\widehat{G},\sigma_U \ast \sigma_V)$ is one-to-one mod 0.
If $U$ and $V$ have simple spectrum then they are strongly disjoint if and only if $U \otimes V$ has simple spectrum, and hance for any two strongly disjoint unitary representations $U$ and $V$ we have $\mc M(U\otimes V) = \mc M(U) \mc M(V)$.
In fact, Lemma~\ref{LemRSS} gives the useful sufficient condition of strong disjointness for unitary representations.

\subsection{Group actions}

Let $\Gamma$ be a locally compact second countable group.
Given a standard non-atomic probability space $(X,\goth B, \mu)$, let $\text{Aux}(X,\mu)$ stand for the group of invertible $\mu$-preserving transformations of $X$.
By an \emph{action} $T$ of $\Gamma$ we mean a continuous group homomorphism $T\colon \Gamma\ni g\mapsto T_g\in\text{Aut}(X,\mu)$.
Denote by $\mc A_\Gamma \subset\text{Aut}(X,\mu)^\Gamma$ the subset of all measure-preserving actions of $\Gamma$ on $(X,\goth B,\mu)$.
Recall that $U_T$ denotes the Koopman representation of $\Gamma$ associated with $T\in\mc A_\Gamma$.
We endow $\mc A_\Gamma$ with the weakest topology which makes continuous the mapping
$$
\mc A_\Gamma\ni T \mapsto U_T \in \mc U_\Gamma(L_0^2(X,\mu)).
$$
It is Polish.
It is easy to verify that a sequence $T^{(n)}$ of $\Gamma$-actions converges to $T$ if and only if $\sup_{g \in K} \mu(T^{(n)}_g A\bigtriangleup T_g A)\to 0$ as $n\to\infty$ for each compact $K\subset \Gamma$  and $A\in \goth B$.
There is a natural action of $\text{Aut}(X,\mu)$ on $\mc A_\Gamma$ by conjugation:
$$
(R\cdot T)_g = R T_g R^{-1} \quad \text{for $R\in\text{Aut}(X,\mu)$, $T\in\mc A_\Gamma$, $g\in \Gamma$,}
$$
and this action is obviously continuous.

If $\mu(X)=\infty$ we define the Polish space $\mc A_\Gamma (X,\mu)$ of all infinite measure preserving $\Gamma$-actions in a similar way.
Notice that for $\mu$ is infinite the Koopman representation associated with $T\in\mc A_\Gamma (X,\mu)$ is considered in the entire space $L^2(X,\mu)$.

\subsection{$(C,F)$-construction} \label{CF_sec}

We now briefly outline the $(C,F)$-construction of measure preserving actions for locally compact groups.
For details see \cite{Da_CF} and references therein.

Let $\Gamma$ be a unimodular locally compact second countable amenable group.
Fix a ($\sigma$-finite) left Haar measure $\lambda$ on it.
Given two subsets $E,F\subset \Gamma$, by $EF$ we mean their algebraic product, i.e. $EF=\{ef \mid e\in E, f\in F\}$.
The set $\{e^{-1} \mid e\in E\}$ is denoted by $E^{-1}$.
If $E$ is a singleton, say $E=\{e\}$, then we will write $eF$ for $EF$.

To define a $(C,F)$-action of $\Gamma$ we need two sequences $(F_n)_{n=0}^\infty$ and $(C_n)_{n=1}^\infty$ of subsets in $\Gamma$ such that the following conditions are satisfied:
\begin{align}
& (F_n)_{n=0}^\infty \text{ is a F{\o}lner sequence in $\Gamma$}, \label{CF1}\\
& C_n \text{ is finite and } \# C_n > 1,  \label{CF2}\\
& F_n C_{n+1} \subset F_{n+1},  \label{CF3}\\
& F_n c \cap F_n c' = \emptyset \text{ for all } c\neq c'\in C_{n+1}.  \label{CF4}
\end{align}
We equip $F_n$ with the measure $(\#C_1\cdots\#C_n)^{-1}\lambda\upharpoonright F_n$ and endow $C_n$ with the equidistributed probability measure.
Let $X_n := F_n \times \prod_{k>n} C_k$ stand for the product of measure spaces.
Define an embedding $X_n \to X_{n+1}$ by setting
$$
(f_n,c_{n+1},c_{n+2},\ldots) \mapsto (f_n c_{n+1},c_{n+2},\ldots).
$$
It is easy to see that this embedding is measure preserving.
Then $X_1\subset X_2\subset\cdots$.
Let $X:=\bigcup_{n=0}^\infty X_n$ denote the inductive limit of the sequence of measure spaces $X_n$ and let $\goth B$ and $\mu$ denote the corresponding Borel $\sigma$-algebra and measure on $X$.
Then $X$ is a standard Borel space with $\mu$ is $\sigma$-finite.
It is finite if
\begin{equation} \label{CFfin}
\prod_{n=1}^\infty \frac{\lambda(F_{n+1})}{\lambda(F_n) \#C_{n+1}} < \infty.
\end{equation}
and infinite if
\begin{equation} \label{CFinf}
\prod_{n=1}^\infty \frac{\lambda(F_{n+1})}{\lambda(F_n) \#C_{n+1}} = \infty.
\end{equation}
If (\ref{CFfin}) is satisfied then we choose (i.e., normalize) $\lambda$ in such a way that $\mu(X)=1$.
Given a Borel subset $A\subset F_n$, we put
$$
[A]_n := \{x\in X \mid x=(f_n,c_{n+1},c_{n+2},\ldots)\in X_n \text{ and } f_n\in A\}
$$
and call this set an \emph{$n$-cylinder}.
It is clear that the $\sigma$-algebra $\goth B$ is generated by the family of all cylinders.

To construct $\mu$-preserving action of $\Gamma$ on $(X,\goth B,\mu)$, fix a filtration $K_1\subset K_2\subset\cdots$ of $\Gamma$ by compact subsets.
Thus $\bigcup_{m=1}^\infty K_m = \Gamma$.
Given $n,m\in\N$, we set
\begin{align*}
L_m^{(n)} &:= \left( \bigcap_{k\in K_m} (k^{-1}F_n)\cap F_n \right) \times \prod_{k>n} C_k \subset X_n \text{ and } \\
R_m^{(n)} &:= \left( \bigcap_{k\in K_m} (k     F_n)\cap F_n \right) \times \prod_{k>n} C_k \subset X_n.
\end{align*}
It is easy to verify that $L_{m+1}^{(n)} \subset L_m^{(n)} \subset L_m^{(n+1)}$ and $R_{m+1}^{(n)} \subset R_m^{(n)} \subset R_m^{(n+1)}$.
We define a Borel mapping $K_m \times L_m^{(n)} \ni (g,x) \mapsto T_{m,g}^{(n)} x \in R_m^{(n)}$ by setting for $x=(f_n,c_{n+1},c_{n+2},\ldots)$,
$$
T_{m,g}^{(n)}(f_n,c_{n+1},c_{n+2},\ldots) := (g f_n,c_{n+1},c_{n+2},\ldots).
$$
Now let $L_m := \bigcup_{n=1}^\infty L_m^{(n)}$ and $R_m := \bigcup_{n=1}^\infty R_m^{(n)}$.
Then a Borel one-to-one mapping $T_{m,g}\colon K_m \times L_m \ni (g,x) \mapsto T_{m,g}x \in R_m$ is well defined by 
$T_{m,g} \upharpoonright L_m^{(n)} = T_{m,g}^{(n)}$ for $g \in K_m$ and $n\geqslant 1$.
It is easy to see that $L_m \supset L_{m+1}$, $R_m \supset R_{m+1}$ and $T_{m,g}\upharpoonright L_{m+1} = T_{m+1,g}$ for all $m$.
It follows from (\ref{CF1}) that $\mu(L_m) = \mu(R_m) = 1$ for all $m\in\N$.
Finally we set $\widehat{X} := \bigcap_{m=1}^\infty L_m \cap \bigcap_{m=1}^\infty R_m$ and define a Borel mapping $T\colon \Gamma\times \widehat{X} \ni (g,x) \mapsto T_g x \in \widehat{X}$ by setting $T_g x := T_{m,g}x$ for some (and hence any) $m$ such that $g\in K_m$.
It is clear that $\mu(\widehat{X})=1$.
Thus we obtain that $T=(T_g)_{g\in \Gamma}$ is a free Borel measure preserving action of $\Gamma$ on a conull subset of a standard Borel space $(X,\goth B,\mu)$.
It is easy to verify that $T$ does not depend on the choice of filtration $(K_m)_{m=1}^\infty$.
$T$ is called the \emph{$(C,F)$-action of $\Gamma$ associated with $(C_{n+1},F_n)_{n\geqslant 0}$}.

We now recall some basic properties of $(X, \goth B, \mu, T)$.
Given Borel subsets $A,B\subset F_n$, we have
\begin{align*}
& [A\cap B]_n = [A]_n \cap [B]_n, [A\cup B]_n = [A]_n \cup [B]_n, \\
& [A]_n = [AC_{n+1}]_{n+1} = \bigsqcup_{c\in C_{n+1}} [Ac]_{n+1}, \\
& T_g[A]_n = [gA]_n \text{ if } gA\subset F_n.
\end{align*}
Note also that the $(C,F)$-construction `respects' Cartesian products.
Namely, the product of two $(C,F)$-actions $(T_g^{(i)})_{g \in G_i}$ associated with $(C_n^{(i)},F_n^{(i)})_n$, $i=1,2$, is the $(C,F)$-action of $G_1\times G_2$ associated with $(C_n^{(1)}\times C_n^{(2)},F_n^{(1)}\times F_n^{(2)})_n$.

\subsection{Poisson suspension} \label{PoisSusp_sec}

Let $(X,\mc B)$ be a standard Borel space and let $\mu$ be an infinite $\sigma$-finite non-atomic measure on $X$.
Fix an increasing sequence of Borel subsets $X_1 \subset X_2 \subset \cdots$ with $\bigcup_{i=1}^\infty X_i = X$ and $\mu(X_i)<\infty$ for each $i$.
A Borel subset is called \emph{bounded} if it is contained in some $X_i$.
Let $\widetilde{X}_i$ denote the space of finite measures on $X_i$.
For each bounded subset $A\subset X_i$, let $N_A$ stand for the map
$$
\widetilde{X}_i \ni \omega \mapsto \omega(A) \in \R.
$$
Denote by $\widetilde{\mc B}_i$ the smallest $\sigma$-algebra on $\widetilde{X}_i$ in which all the maps $N_A$, $A\in \mc B\cap X_i$, are measurable.
It is well known that $(\widetilde{X}_i,\widetilde{\mc B}_i)$ is a standard Borel space.
Denote by $(\widetilde{X},\widetilde{\mc B})$ the projective limit of the sequence
$$
(\widetilde{X}_1,\widetilde{\mc B}_1) \leftarrow (\widetilde{X}_2,\widetilde{\mc B}_2) \leftarrow \cdots,
$$
where the arrows denote the (Borel) natural restriction maps.
Then $(\widetilde{X},\widetilde{\mc B})$ is a standard Borel space.
To put it in other way, $\widetilde{X}$ is the space of measures on $X$ which are $\sigma$-finite along $(X_i)_{i>0}$.
Then there is a unique probability measure $\widetilde{\mu}$ on $(\widetilde{X},\widetilde{\mc B})$ such that
\begin{enumerate}
\item
$N_A$ maps $\widetilde{\mu}$ to the Poisson distribution with parameter $\mu(A)$, i.e.
$$
\widetilde{\mu}(\{\omega \mid N_A(\omega)=j\}) = \frac{\mu(A)^j \exp(-\mu(A))}{j!}
$$
for all bounded $A\subset X$ and integer $j\geqslant 0$ and
\item
if $A$ and $B$ are disjoint bounded subsets of $X$ then the random variables $N_A$ and $N_B$ on $(\widetilde{X}, \widetilde{\mc B}, \widetilde{\mu})$ are independent.
\end{enumerate}

Let $G$ be a locally compact second countable group and let $T$ be a $\mu$-preserving action of $G$ on $X$ such that $T_g$ preserves the subclass of bounded subsets for each $g\in G$.
Then $T$ induces a $\widetilde{\mu}$-preserving action $\widetilde{T}$ of $G$ on $\widetilde{X}$ by the formula $\widetilde{T}_g\omega := \omega\circ T_{-g}$.
We recall that the dynamical system $(\widetilde{X},\widetilde{\mc B},\widetilde{\mu},\widetilde{T})$ is called the \emph{Poisson suspension} of $(X,\mc B,\mu,T)$ (see \cite{CFS}, \cite{Roy2} for the case $G=\Z$).

The well known Fock representation of $L^2(\widetilde{X},\widetilde{\mu})$ gives an isomorphism
$$
L^2(\widetilde{X},\widetilde{\mu}) \simeq \bigoplus_{n=0}^\infty L^2(X,\mu)^{\odot n},
$$
where $L^2(X,\mu)^{\odot n}$ is the $n$th symmetric tensor power of $L^2(X,\mu)$, with $L^2(X,\mu)^{\odot 0} = \CC$.
The Koopman representation $U_{\widetilde{T}}\oplus P_0$ (considered on $L^2(\widetilde{X},\widetilde{\mu})$) is unitarily equivalent to the exponential of $U_T$:
$$
U_{\widetilde{T}}\oplus P_0 \simeq \exp U_T = \bigoplus_{n=0}^\infty U_T^{\odot n},
$$
where $P_0$ is the orthogonal projection on $\CC \subset L^2(\widetilde{X},\widetilde{\mu})$ and $U_T^{\odot n}$ is the $n$th symmetric tensor power of $U_T$ \cite{Ne}.
Recall that since $\mu$ is infinite, we consider $U_T$ in the entire space $L^2(X,\mu)$.
It follows, in particular, that the mapping $\mc A_\Gamma (X,\mu) \ni T \mapsto \widetilde{T} \in \mc A_\Gamma(\widetilde{X},\widetilde{\mu})$ is continuous.
$\widetilde{T}$ is rigid (for the sequence $g_n$) if and only if $T$ is rigid (for the sequence $g_n$).
If $T$ has no invariant subsets of finite positive measure then $\widetilde{T}$ is weakly mixing \cite{Roy2}.

\section{$\R^m$-actions} \label{Rp_sec}

In this section we prove Theorem~\ref{MainTh} in the case when $G=\R^m$.

For given $p > 1$, let $A \colon \Z^p \to \Z^p$ denote a `cyclic' group automorphism: $$A(x_1, x_2, \ldots, x_p) = (x_p, x_1, \ldots, x_{p-1}).$$
Following~\cite{DS}, denote by $\Gamma$ the semidirect product\footnote{By $\Z(p)$ we denote a cyclic group of order $p$, i.e. $\Z(p)=\Z/p\Z$.}
$$
\Gamma := G \times \Z^p \rtimes_A \Z(p)
$$
with the multiplication law as follows:
$$
(g,x,n)(h,y,k) := (g+h, x+A^n y, n+k),
$$
$g,h \in G$, $x,y \in \Z^p$, $n,k \in \Z(p)$.
We will identify $G$ with the subgroup $\{(g,0,0) \mid g\in G\} \subset \Gamma$.
Let $\mc E_\Gamma \subset \mc A_\Gamma$ stand for the subset of all free ergodic $\Gamma$-actions.
$\mc E_\Gamma$ is $G_\delta$ subset in $\mc A_\Gamma$ and hence it is Polish group with the induced topology \cite{DS}.
To prove Theorem~\ref{MainTh} we will use `generic' argument and the following facts will be needed.

\begin{lemma}[{\cite[Theorem~2.8]{DS}}] \label{LemHS}
For a generic action $T\in \mc E_\Gamma$ the action $T\upharpoonright G$ is weakly mixing and $\mc M(T\upharpoonright G) = \{p\}$.
\end{lemma}

\begin{lemma}[{\cite[Lemma~2.4]{DS}}] \label{dense}
The $\text{Aut}(X,\mu)$-orbit of any action $T\in \mc E_\Gamma$ is dense in $\mc E_\Gamma$.
\end{lemma}

We will apply Lemma~\ref{dense} to show that the set of $\Gamma$-actions with certain properties is dense in $\mc E_\Gamma$.
However to apply this lemma we will need at least one action in this set.
This single action is constructed explicitly in Lemma~\ref{LemRigid}.

\begin{lemma} \label{LemRigid}
For any sequence $(g_k)_{k=1}^\infty \subset G$, $g_k \to \infty$, there exists a $(C,F)$-action $T \in \mc E_\Gamma$ such that $U_T(g_{k_n}) \to I$ for some subsequence $(g_{k_n})_{n=1}^\infty$ of $(g_k)_{k=1}^\infty$.
\end{lemma}

\begin{proof}
To construct $(C,F)$-action we shall determine a sequence $(C_{n+1},F_n)_{n=0}^\infty$.
This will be done inductively.
Let $F_n = F'_n  \times F''_n $ and $C_n = C'_n  \times C''_n $, where $F'_n , C'_n  \subset G$, $F''_n , C''_n  \subset \Z^p \rtimes \Z(p)$.

First, we claim that the sets $F'_n , C'_n  \subset G = \R^m$ and a subsequence $(g_{k_n})_{n=1}^\infty$ can be chosen in such a way that
$$
\lim_{n\to\infty} \frac{\#(C'_n  \cap (C'_n  - g_{k_n}))}{\# C'_n } = 1.
$$
Let us show this.
Select a subsequence $(g_{k_n})_{n=1}^\infty$ of $(g_k)_{k=1}^\infty$ such that $\dfrac{g_{k_n}}{|g_{k_n}|}$ converges (to some point of the unit sphere) as $n \to \infty$.
From now on we will write $g_n$ instead of $g_{k_n}$ for short.
Let $g_n=(g_n^{(1)},\ldots,g_n^{(m)}) \in R^m$.
Without loss of generality we may assume that $g_n^{(i)} > 0$, $i=1,\ldots,m$, and $g_n^{(1)} \to \infty$.
In the other cases the proof is similar.
Fix a sequence of positive integers $\alpha_n$ with $\sum_{n=1}^\infty \alpha_n < \infty$.
By replacing $(g_n)_{n=1}^\infty$ with its subsequence if necessary, we may assume that
$$
\frac{g_{n+1}^{(1)}}{g_n^{(1)}} > \frac{1}{\alpha_n} + 1.
$$
We will construct $C'_n $ and $F'_n $ inductively.
Choose $C'_n $ and $F'_n $ arbitrarily.
Now suppose that we already have $C'_{n-1}$ and
$F'_{n-1} = (-a_{n-1}^{(1)},a_{n-1}^{(1)}) \times \cdots \times (-a_{n-1}^{(m)},a_{n-1}^{(m)})$,
where $a_{n-1}^{(i)} > 0$ and $a_{n-1}^{(1)} = \dfrac{g_n^{(1)}}{2}$.
Our purpose is to define $C'_n $ and $F'_n $.
Set 
$$
h_n := \left\lfloor \frac{g_{n+1}^{(1)}-g_n^{(1)}}{2g_n^{(1)}} \right\rfloor > \frac{1}{\alpha_n}-\frac{1}{2}.
$$
In particular, $(2 h_n + 1) g_n^{(1)} < g_{n+1}^{(1)} < (2 h_n + 1) g_n^{(1)} + 2 g_n^{(1)}$ and $2h_n+1 > \dfrac{2}{\alpha_n}$.
Select integers $w_n^{(2)},\ldots,w_n^{(m)} > 0$ in such a way that
$$
\frac{h_n g_n^{(i)}}{a_{n-1}^{(1)}(2w_n^{(i)}+1)} < \alpha_n.
$$
We set
\begin{align*}
A_n & := \{ (0, 2l^{(2)}a_{n-1}^{(2)}, \ldots, 2l^{(m)}a_{n-1}^{(m)}) \mid l^{(i)} \in \Z, -w_n^{(i)} \leq l^{(i)} \leq w_n^{(i)} \} \subset \R^m, \\
C'_n & := \bigsqcup_{k=-h_n}^{h_n} (A_n + k g_n).
\end{align*}
Let also $a_n^{(1)} := \dfrac{g_{n+1}^{(1)}}{2}$ and $a_n^{(i)} := (2 w_n^{(i)} + 1) a_{n-1}^{(i)} + h_n g_n^{(i)}$, $i=2,\ldots,m$.
Set $F'_n := (-a_n^{(1)},a_n^{(1)}) \times \cdots \times (-a_n^{(m)},a_n^{(m)})$.
Then, by construction,
\begin{multline*}
\frac{\lambda(F'_n)}{\lambda(F'_{n-1}) \#C'_n} = \frac{g_{n+1}^{(1)}}{g_n^{(1)}(2 h_n + 1)} \prod_{i=2}^m \frac{2a_n^{(i)}}{2a_{n-1}^{(i)}(2w_n^{(i)} + 1)} < \\
< \left( 1 + \frac{2}{2h_n + 1} \right) \prod_{i=2}^m \left( 1 + \frac{h_n g_n^{(i)}}{2w_n^{(i)} + 1} \right) < (1+\alpha_n)^m,
\end{multline*}
Thus the conditions~(\ref{CF1})-(\ref{CFfin}) hold for $(F'_n ,C'_n )_n$.
It also follows from the definition of $C'_n $ that
$$
\frac{\# (C'_n  \cap (C'_n  - g_n))}{\# C'_n } = \frac{2h_n}{2h_n+1} \to 1.
$$

Secondly, let 
$C''_n $ and $F''_n $ be any subsets of $\Z^p\rtimes \Z(p)$ satisfying~(\ref{CF1})-(\ref{CFfin}).
For instance, 
set
\begin{align*}
F''_n  &:= \left\{ -\tfrac{3^n-1}{2} ,\ldots, \tfrac{3^n-1}{2} \right\}^p \times \Z(p) \subset \Z^p\rtimes \Z(p), \\
C''_n  &:= \left\{-3^{n-1},0, 3^{n-1}\right\}^p \times \{0\} \subset \Z^p\rtimes \Z(p).
\end{align*}

Let $T$ be $(C,F)$-action associated with $(C_n,F_n)_n = (C'_n\times C''_n, F'_n\times F''_n)_n$.
As was mentioned in Section~\ref{CF_sec}, $T$ is then the product of two $(C,F)$-actions $T^{(1)}=(T^{(1)}_g)_{g \in G}$ and $T^{(2)}=(T^{(2)}_z)_{z \in \Z^p\rtimes \Z(p)}$ associated with $(C'_n ,F'_n )_n$ and $(C''_n ,F''_n )_n$ respectively.
Since $g_n \in G$, we have
$$
\lim_{n\to\infty} \frac{\# (C_n \cap g_n^{-1} C_n)}{\# C_n} = 1.
$$

We claim that $\lim_{n \to \infty} \mu(T_{g_n}A \bigtriangleup A) = 0$ for any $A \in \goth B$.
It suffices to consider the  cylinders $[A]_n$, $A \subset F_n$.
Fix arbitrary $\varepsilon > 0$ and select $n$ such that
\begin{equation} \label{eps}
\frac{\#(C_n \setminus g_n^{-1} C_n)}{\# C_n} < \varepsilon.
\end{equation}
Let $A \subset F_{n-1}$.
Notice that $g_n$ commutes with all the elements of $\Gamma$.
Thus we have
$$
[A]_{n-1} = \bigsqcup_{c \in C_n} [Ac]_n = A_1 \sqcup \bigsqcup_{c \in C_n\cap g_n C_n} [Ac]_n = A_2 \sqcup \bigsqcup_{c \in C_n\cap g_n^{-1} C_n} [Ac]_n,
$$
where $A_1 := \bigsqcup_{c \in C_n\setminus g_n C_n} [Ac]_n$, $A_2 := \bigsqcup_{c \in C_n\setminus g_n^{-1} C_n} [Ac]_n$ and $\mu(A_i) < \varepsilon$ by~(\ref{eps}).
On the other hand,
\begin{align*}
T_{g_n} [A]_{n-1} &= T_{g_n} A_2 \sqcup \bigsqcup_{c \in C_n\cap g_n^{-1} C_n} T_{g_n} [Ac]_n \\
&= T_{g_n} A_2 \sqcup \bigsqcup_{c \in C_n\cap g_n^{-1} C_n} [g_n Ac]_n \\
&= T_{g_n} A_2 \sqcup \bigsqcup_{c \in C_n\cap g_n C_n} [Ac]_n.
\end{align*}
Hence $T_{g_n} [A]_{n-1} \bigtriangleup [A]_{n-1} \subset A_1 \cup T_{g_n} A_2$ and $\mu( T_{g_n} [A]_{n-1} \bigtriangleup [A]_{n-1} ) < 2\varepsilon$.
The claim is proven.
It follows that $U_T(g_n) \to I$ as $n\to\infty$.

Since any $(C,F)$-action is free and ergodic, $T\in\mc E_\Gamma$.
\end{proof}


As was mentioned above, to prove the main result we will apply the Baire category theorem, so the following lemma will be useful.

\begin{lemma} \label{limI}
Given a sequence $g_n \to \infty$ in $G$, the following subsets are residual in $\mc E_\Gamma$:
\begin{align*}
\mc I &:= \{ T \in \mc E_\Gamma \mid I \text{ is a limit point of } \{U_T(g_n)\}_{n=1}^\infty \} \text{ and } \\
\mc O &:= \{ T \in \mc E_\Gamma \mid 0 \text{ is a limit point of } \{U_T(g_n)\}_{n=1}^\infty \}.
\end{align*}
\end{lemma}

\begin{proof}
It follows from Lemma~\ref{Pol} that $\mc I$ and $\mc O$ are both $G_\delta$ in $\mc E_\Gamma$.
Notice also that $\mc I$ and $\mc O$ are both $\text{Aut}(X,\mu)$-invariant.
Therefore in view of Lemma~\ref{dense} it remains to show that $\mc I$ and $\mc O$ contain at least one free ergodic action.
$\mc I$ is non-empty by Lemma~\ref{LemRigid}.
Consider an action of $\Gamma$ on itself by translations.
This action preserves the ($\sigma$-finite, infinite) Haar measure.
The corresponding Poisson suspension (see Section~\ref{PoisSusp_sec}) of this action is a probability preserving free $\Gamma$-action and it belongs to $\mc O$ (see \cite{OW}).
\end{proof}


Lemma~\ref{MailLemR} will be the main ingredient in the proof of Theorem~\ref{MainTh}.
In general, its proof goes along the lines developed in \cite{Ry} for $\Z$-actions.

\begin{lemma} \label{MailLemR}
Given a rigid weakly mixing $S\in \mc A_G$ and $p > 0$, there exists a weakly mixing $T\in \mc A_G$ such that $S \times T$ is rigid, weakly mixing and $\mc M(S \times T) = \mc M(S) \diamond \{p\}$.

Moreover, if $(r_n)_{n=1}^\infty$ and $(g_n)_{n=1}^\infty$ are sequences in $G$ such that $U_S (r_n) \to I$ and $U_S(g_n)\to 0$, then $U_{S\times T} (r'_n) \to I$ and $U_{S\times T}(g'_n)\to 0$ as $n\to\infty$ for some subsequences $(r'_n)_{k=1}^\infty$ and $(g'_n)_{k=1}^\infty$.
\end{lemma}

\begin{proof}
Fix sequences $(r_n)_{n=1}^\infty$ and $(g_n)_{n=1}^\infty$ in $G$ such that
\begin{align}
U_S (r_n) &\to I \text{ and} \label{SrnI} \\
U_S (g_n) &\to 0.  \label{Srn0}
\end{align}
Let also $\langle d_1,\ldots,d_{2m} \rangle$ be dense in $G$.
Let $\Gamma := G \times \Z^p \rtimes_A \Z(p)$ stand for the auxiliary non-Abelian group defined above.
We claim that for a generic $\widetilde{T} \in \mc E_\Gamma$, $G$-action $T:=\widetilde{T}\upharpoonright G$ satisfies the following properties:
\begin{enumerate}
\item \label{wm}
$T$ is weakly mixing,
\item \label{hs}
$\mc M(T) = \{p\}$,
\item \label{rU}
$0$, $I$ and $U_T(d_1),\ldots,U_T(d_{2m})$ are limit points of the set $\{U_T(r_n)\}_{n\in\N}$,
\item \label{gI}
$0$ and $I$ are limit points of $\{U_T(g_n)\}_{n\in\N}$.
\end{enumerate}
The properties (\ref{wm})--(\ref{hs}) are generic by Lemma~\ref{LemHS}.
Since $U_T(d)$ is a limit point of $\{U_T(r_n)\}_{n=1}^\infty$ if and only if $I$ is a limit point of $\{U_T(r_n-d)\}_{n=1}^\infty$, Lemma~\ref{limI} implies (\ref{rU})--(\ref{gI}) for a generic $\widetilde{T} \in \mc E_\Gamma$.
Hence there is an action satisfying all of these conditions.

Now let us show that $T$ is the required action.
Lemma~\ref{LemRSS}, in view of (\ref{SrnI}) and (\ref{rU}), implies that $\mc M(U_S \otimes U_T) = p \mc M(U_S)$.
Since the Koopman representation is considered on the space $L^2(X,\mu)\ominus \CC$, we have
\begin{equation} \label{tensum}
U_{S \times T} = (1 \otimes U_T)\oplus(U_S \otimes U_T)\oplus(U_S \otimes 1),
\end{equation}
where $1$ denotes the identity operator on $\CC$.
If $1 \otimes U_T$, $U_S \otimes U_T$, $U_S \otimes 1$ are pairwise spectrally disjoint then
$$
\mc M(S \times T) = \{p\} \cup p \mc M(S) \cup \mc M(S) = \mc M(S) \diamond \{p\}.
$$
Apply (\ref{rU}) and (\ref{gI}) and fix a subsequence $(r''_n)_{n=1}^\infty$ of $(r_n)_{n=1}^\infty$ and a subsequence $(g''_n)_{n=1}^\infty$ of $(g_n)_{n=1}^\infty$ such that $U_T(r''_n)\to 0$ and $U_T(g''_n)\to I$ as $n\to\infty$.
The spectrally disjointness for each pair of terms from (\ref{tensum}) follows frow Lemma~\ref{LemDis}, since
$$
(U_S \otimes 1)(r''_n) \to I, \quad (U_S \otimes U_T)(r''_n) \to 0,
$$
$$
(1 \otimes U_T)(g''_n) \to I, \quad (U_S \otimes U_T)(g''_n) \to 0,
$$
$$
(1 \otimes U_T)(g''_n) \to I, \quad (U_S \otimes 1)(g''_n) \to 0.
$$

It is clear that $S \times T$ is weakly mixing.
By (\ref{rU}) and (\ref{gI}) there are subsequences $(r'_n)_{n=1}^\infty$ and $(g'_n)_{n=1}^\infty$ of $(r_n)_{n=1}^\infty$ and $(g_n)_{n=1}^\infty$ such that $U_T(r'_n)\to I$ and $U_T(g'_n)\to 0$.
Hence $U_{S\times T} (r'_n) \to I$ and $U_{S\times T}(g'_n)\to 0$.
\end{proof}



\begin{proof}[Proof of Theorem~\ref{MainTh} for $G=\R^m$]
Consider the auxiliary group $\Gamma_1 := G \times \Z^{p_1} \rtimes \Z(p_1)$ defined above.
Let $\widetilde{T}_1 \in\mc E_{\Gamma_1}$ be such that $T_1 := \widetilde{T}_1 \upharpoonright G$ is weakly mixing, $\mc M(T_1) = \{p_1\}$ and $U_{T_1}(r_{n,1})\to I$, $U_{T_1}(g_{n,1})\to 0$ as $n \to \infty$, where $(r_{n,1})_{n=1}^\infty$, $(g_{n,1})_{n=1}^\infty$ are some sequences in $G$.
Since all these properties are generic for the actions from $\mc E_{\Gamma_1}$ by Lemmata~\ref{LemHS} and \ref{limI}, there is an action $\widetilde{T}_1$ possessing all of them.

Now we apply Lemma~\ref{MailLemR} and choose a weakly mixing $T_2 \in \mc A_G$ such that $\mc M(T_1 \times T_2) = \{p_1\}\diamond\{p_2\}$ and $U_{T_1 \times T_2}(r_{n,2})\to I$, $U_{T_1 \times T_2}(g_{n,2})\to 0$ as $n \to \infty$, where $(r_{n,2})_{n=1}^\infty$ and $(g_{n,2})_{n=1}^\infty$ are subsequences of $(r_{n,1})_{n=1}^\infty$ and $(g_{n,1})_{n=1}^\infty$ respectively.


By induction, given a weakly mixing $G$-action $T_1\times\cdots\times T_{k-1}$ with
$$\mc M(T_1 \times\cdots\times T_{k-1}) = \{p_1\}\diamond\cdots\diamond \{p_{k-1}\},$$
$$U_{T_1\times\cdots\times T_{k-1}}(r_{n,k-1})\to I,$$
$$U_{T_1\times\cdots\times T_{k-1}}(g_{n,k-1})\to 0,$$
by Lemma~\ref{MailLemR} there exists a weakly mixing $T_k \in \mc A_G$ such that
\begin{equation} \label{MTR1}
\mc M(T_1 \times\cdots\times T_k) = \{p_1\}\diamond\cdots\diamond \{p_k\},
\end{equation}
\begin{equation} \label{MTR2}
U_{T_1\times\cdots\times T_k}(r_{n,k})\to I,
\end{equation}
\begin{equation} \label{MTR3}
U_{T_1\times\cdots\times T_k}(g_{n,k})\to 0,
\end{equation}
where $(r_{n,k})_{n=1}^\infty$ and $(g_{n,k})_{n=1}^\infty$ are subsequences of $(r_{n,k-1})_{n=1}^\infty$ and $(g_{n,k-1})_{n=1}^\infty$ respectively.
This proves the theorem in the case when the sequence $p_1,p_2,\ldots$ is finite.
Otherwise we obtain an infinite sequence of weakly mixing $G$-actions $T_k$ satisfying (\ref{MTR1})--(\ref{MTR3}).
It is clear that the product $T:=T_1\times T_2 \times\cdots$ is weakly mixing and $\mathcal M(T) = \{p_1\}\diamond\{p_2\}\diamond\cdots$.
\end{proof}


The following simple lemma (that was stated in \cite{DL} without proof) shows how to extend the result of Theorem~\ref{MainTh} from $\R$ to any torsion free discrete countable Abelian group (Corollary~\ref{CorTorFr}).

\begin{lemma} \label{emb}
Let $G$ and $H$ be locally compact second countable Abelian groups and let $\varphi\colon G\to H$ be a continuous one-to-one homomorphism with $\overline{\varphi(G)}=H$.
Given an $H$-action $T=(T_h)_{h\in H}$, the composition $T\circ\varphi = (T_{\varphi(g)})_{g\in G}$ is $G$-action with $\mc M(T\circ\varphi) = \mc M(T)$. 
\end{lemma}

\begin{proof}
Let $\sigma$ be a measure of maximal spectral type and $m\colon \widehat{H}\to \N\cup\{\infty\}$ be the spectral multiplicity function of $U_T$:
\begin{equation}\label{L2decomp}
L^2_0(X,\mu) = \int_{\widehat{H}}^\oplus \mc H_\chi d\sigma(\chi) 
\quad \text{and}\quad
U_T(h)f(\chi) = \chi(h) f(\chi),\; h\in H, 
\end{equation}
for each $f\colon \widehat{H} \ni \chi \mapsto f(\chi) \in \mc H_\chi$ with $\int_{\widehat{H}} \|f(\chi)\|^2 d\sigma(\chi) < \infty$, $\dim \mc H_\chi = m(\chi)$.
Let $\widehat{\varphi}\colon \widehat{H}\to \widehat{G}$ stand for the dual to $\varphi$ homomorphism and $\widehat{\sigma}:=\sigma\circ\widehat{\varphi}^{-1}$ be the image of $\sigma$ with respect to $\widehat{\varphi}$.
Clearly, $\widehat{\sigma}(\widehat{\varphi}(\widehat{H}))=1$.
Let $\sigma=\int_{\widehat{G}} \sigma_\omega d\widehat{\sigma}(\omega)$ denote the disintegration of $\sigma$ relative to $\widehat{\varphi}$.
Then we derive from (\ref{L2decomp}) that
$$
L^2_0(X,\mu) = \int_{\widehat{G}}^\oplus \mc H'_\omega d\widehat{\sigma}(\omega)
 = \int_{\widehat{\varphi}(\widehat{H})}^\oplus \mc H'_\omega d\widehat{\sigma}(\omega),
$$
where $\mc H'_\omega:=\int_{\widehat{H}}^\oplus \mc H_\chi d\sigma_\omega(\chi)$.
Let $l(\omega) := \dim\mc H_\omega$, $\omega\in\widehat{G}$.
Then
$$
l(\omega) =
\begin{cases}
\sum_{\sigma_\omega(\chi)>0} m(\chi),& \text{if $\sigma_\omega$ is purely atomic,}\\
\infty, & \text{otherwise.}
\end{cases}
$$
Since $\overline{\varphi(G)}=H$, $\widehat{\varphi}$ is one-to-one and hance $\mc H'_{\widehat{\varphi}(\chi)} = \mc H_\chi$  for any $\chi\in\widehat{H}$.
In particular, $l(\widehat{\varphi}(\chi)) = m(\chi)$, $\chi\in\widehat{H}$.
It follows from (\ref{L2decomp}) that for any $\omega = \widehat{\varphi}(\chi) \in \widehat{G}$
\begin{multline*}
U_{T\circ\varphi}(g)f(\omega) =
U_T(\varphi(g))f(\widehat{\varphi}(\chi)) =
U_T(\varphi(g))(f\circ\widehat{\varphi})(\chi) = \\
= \chi(\varphi(g))(f\circ\widehat{\varphi})(\chi) =
(\widehat{\varphi}(\chi))(g)f(\widehat{\varphi}(\chi)) =
\omega(g)f(\omega).
\end{multline*}
This means that $\widehat{\sigma}$ is a measure of maximal spectral type and $l$ is the spectral multiplicity function of $U_{T\circ\varphi}$.
Hence $\mc M(T\circ\varphi) = \mc M(T)$.
\end{proof}

\begin{corollary} \label{CorTorFr}
Let $G$ be a torsion free discrete countable Abelian group.
Given a sequence of positive integers $p_1,p_2,\ldots$, there exists a weakly mixing probability preserving $G$-action $S$ such that $\mc M(S) = \{p_1\}\diamond\{p_2\}\diamond\cdots$.
\end{corollary}

\begin{proof}
In the case when $G=\Z$ see \cite{Ry} or Section~\ref{Disc_sec}.
Consider the case when $G\neq\Z$.
In view of Lemma~\ref{emb} it suffices to show that there is an embedding $\varphi\colon G\to\R$ such that $\overline{\varphi(G)}=\R$.
Indeed, $G$ can be embedded into $\Q^\N$ (see \cite{HR}).
In turn, the latter group obviously embeds into $\R$.
It remains to note that if an infinite subgroup of $\R$ is not isomorphic to $\Z$ then it is dense in $\R$.

By Theorem~\ref{MainTh} for $G=\R$, there is a weakly mixing $\R$-action $T$ such that $\mc M(T) = \{p_1\}\diamond\{p_2\}\diamond\cdots$.
Then by Lemma~\ref{emb} the composition $T\circ\varphi = (T_{\varphi(g)})_{g\in G}$ is a weakly mixing $G$-action with $\mc M(T\circ\varphi) = \mc M(T) = \{p_1\}\diamond\{p_2\}\diamond\cdots$.
\end{proof}

\section{Discrete countable Abelian group actions} \label{Disc_sec}

In this section we prove Theorem~\ref{MainTh} in the case when $G$ is an infinite discrete countable Abelian group.

As in the previous section, given a countable discrete Abelian group $J$ and $p>1$, we denote by $\Gamma$ the semidirect product $\Gamma := G \times J^p \rtimes_A \Z(p)$, where $A\colon J^p \to J^p$ is the same (as in Section\ref{Rp_sec}) `cyclic' group automorphism.
From now on we will identify $G$ with the corresponding subgroup in $\Gamma$.

\begin{lemma}[{\cite[Theorem~1.7]{DS}}] \label{DSdisc}
Given $G$ and $p>1$, there is $J$ such that for a generic action $T$ from $\mc A_\Gamma$ the action $T\upharpoonright G$ is weakly mixing and $\mc M(T\upharpoonright G) = \{p\}$.
\end{lemma}
Notice that we can choose $J$ to be either $\Z$ or $\Z(q)^{\oplus\N}$, $q>1$ \cite[Section~1]{DS}.

Let $(g_n)_{n=1}^\infty$ be a sequence in $G$.
We will say that $(g_n)_{n=1}^\infty$ is \emph{`good'} if $g_n \to \infty$ and one of the following is satisfied:
\begin{enumerate}
\item \label{good1}
there is $g_0\in G$ such that $g_n\in\langle g_0\rangle$ for each $n$ (it follows that $g_0$ has infinite order),
\item \label{good2}
each $g_n$ is an element of finite order and orders of $g_n$ are unbounded,
\item \label{good3}
orders of $g_n$ are bounded from above and $g_n$ are independent\footnote{that is, the subgroups $\langle g_n \rangle$ are independent.}.
\end{enumerate}
It is clear that $G$ always contains a `good' sequence.
Notice also that any subsequence of a `good' sequence is also `good'.
We need this notion to be able to apply $(C,F)$-construction in the proof of Lemma~\ref{LemDIlim} which is the analog of Lemma~\ref{LemRigid}.

\begin{lemma} \label{LemDIlim}
Let $(g_k)_{k=1}^\infty$ be a `good' sequence in $G$.
For any $d\in G$ there exists a free action $S\in\mc A_\Gamma$ such that $U_S(g_{k_n}) \to U_S(d)$ for some subsequence $(g_{k_n})_{n=1}^\infty$ of $(g_k)_{k=1}^\infty$.
\end{lemma}

\begin{proof}
Fix $d\in G$.
First, we claim that there is an \textit{infinite} measure preserving action $T$ of $\Gamma$ and subsequence $(g_{k_n})_{n=1}^\infty$ of $(g_k)_{k=1}^\infty$ such that $U_T(g_{k_n}) \to U_S(d)$.
Recall that for $\mu$ infinite, we consider $U_T$ in the entire space $L^2(X,\mu)$.
We will construct $T$ in the form $T=T^{(1)}\times T^{(2)}$, where $T^{(1)}$ and $T^{(2)}$ are $(C,F)$-actions of $G$ and $J^p \rtimes \Z(p)$ respectively.

To construct $T^{(1)}$ we will select subsets $C_n,F_n\subset G$ and a subsequence $(g_{k_n})_{n=1}^\infty$ of $(g_k)_{k=1}^\infty$ in such a way that
\begin{equation}
\lim_{n\to\infty} \frac{\#(C_n \cap (C_n - (g_{k_n}-d)))}{\# C_n} = 1. \label{ArithProgrCond}
\end{equation}
Then, arguing as in the proof of Lemma~\ref{LemRigid}, the reader can easily deduce that $$\lim_{n \to \infty} \mu(T_{g_{k_n}-d}A \bigtriangleup A) = 0$$ for any $A \in \goth B$, and hence $U_T(g_{k_n}-d) \to I$ as $n \to \infty$.

Thus our aim is to select $C_n$, $F_n$ and $k_n$ satisfying (\ref{CF1})--(\ref{CF4}), (\ref{CFinf}) and (\ref{ArithProgrCond}).
This will be done inductively.
Fix an increasing sequence of positive integers $h_n$.
Suppose that we already have $F_{n-1}$ and $k_{n-1}$.
To satisfy (\ref{ArithProgrCond}) we want $C_n$ to be an arithmetic progression with common difference $g_{k_n}-d$ long enough.
We also need $C_n$ to be independent of $F_{n-1}$.
Consider separately three possible cases for $(g_k)_{k=1}^\infty$.

(i)
There is $g_0\in G$ such that $g_k = m_k g_0$, $m_k\in\Z$, $k\in \N$.
Without loss of generality we may assume that $m_k > 0$ and $m_{k+1} > m_k$, $k\in \N$.
Then let $k_n := \max\{ k \mid g_k \in F_{n-1} - F_{n-1} \} + 1$.
Clearly, $l g_{k_n}\notin F_{n-1} - F_{n-1}$ for any $l>0$ and hance $l(g_{k_n}-d) +F_{n-1} \neq l'(g_{k_n}-d) +F_{n-1}$ for $l\neq l'$.

(ii)
Each $g_k$ is an element of finite order and orders of $g_k$ are not bounded.
Without loss of generality we may assume that $\#\{k\mid \ord g_k < N \} < \infty$ for each $N>0$.
Given $0\neq f\in F_{n-1} - F_{n-1}$ and $0 < l \leqslant h_n$, let $D_{n,l}^{f} := \{ k>k_{n-1} \mid l(g_k-d) = f \}$.
We claim that each $D_{n,l}^{f}$ is finite.
Indeed, if $l(g_k-d) = f$ for some $k$ then for any $k'$ with $\ord g_{k'} > l \ord g_k$ we have $\ord (g_k - g_{k'})>l$ and hance $l(g_{k'}-d) \neq l(g_k-d) = f$.
Since there is only finite set of $k'$ with $\ord g_{k'}\leqslant l\ord g_k$,
$D_{n,l}^{f}$ is finite and we can choose $k_n > k_{n-1}$ such that
$k_n \notin D_{n,l}^{f}$ for $0\neq f\in F_{n-1} - F_{n-1}$, $0 < l \leqslant h_n$.
Then $l g_{k_n}\notin F_{n-1} - F_{n-1}$, $0 < l \leqslant h_n$.
In particular, $l(g_{k_n}-d) +F_{n-1} \neq l'(g_{k_n}-d) +F_{n-1}$ for $0 \leqslant l < l' \leqslant h_n$.

(iii)
Orders of $g_k$ are bounded from above and $g_k$ are independent.
In this case for any $0\neq f\in F_{n-1} - F_{n-1}$ and $l>0$ there is at most one $k$ with $l g_k = f$.
Hance we can select $k_n > k_{n-1}$ in such a way that $l g_{k_n} \notin F_{n-1}-F_{n-1}$ whenever $l g_{k_n} \neq 0$.

In each of these three cases we set 
$$
C_n :=
\begin{cases}
\{0,(g_{k_n}-d),2(g_{k_n}-d),\ldots,h_n (g_{k_n}-d)\}, \text{ if $\ord (g_{k_n}-d) > h_n$,}\\
\langle g_{k_n}-d \rangle, \text{ otherwise.}
\end{cases}
$$
It follows that $C_n$ and $F_{n-1}$ are independent.
Since
$$
\frac{\#(C_n \cap (C_n - (g_{k_n}-d)))}{\# C_n} \leqslant \frac{h_n}{h_n+1},
$$
$C_n$ satisfy (\ref{ArithProgrCond}).
Let $F_n\subset G$ be any subset satisfying (\ref{CF1}), (\ref{CF4}) and (\ref{CFinf}).
Let $T^{(1)}$ be $(C,F)$-action associated with $(C_n,F_n)_n$.

$T^{(2)}$ may be any $(C,F)$-action of $J^p \rtimes \Z(p)$.
In view of the structure of $J$ which is either $\Z$ or $\Z(q)^{\oplus\N}$, $q>1$, such an action can be easily constructed.
For instance, set
\begin{align*}
F'_n &:= \left\{ -\tfrac{3^n-1}{2} ,\ldots, \tfrac{3^n-1}{2} \right\}^p \times \Z(p) \subset J^p\rtimes \Z(p), \\
C'_n &:= \left\{-3^{n-1},0, 3^{n-1}\right\}^p \times \{0\} \subset J^p\rtimes \Z(p),
\end{align*}
if $J=\Z$, and
\begin{align*}
F'_n &:= \bigl( \underbrace{\Z(q) \oplus \cdots \oplus \Z(q)}_n \oplus \{0\} \oplus \cdots \bigr)^p \times \Z(p) \subset J^p\rtimes \Z(p), \\
C'_n &:= \bigl( \underbrace{\{0\} \oplus \cdots \oplus \{0\}}_{n-1} \oplus \Z(p) \oplus \{0\} \oplus \cdots \bigr)^p \times \{0\} \subset J^p\rtimes \Z(p)
\end{align*}
if $J = \bigoplus_{n=1}^\infty \Z(q)$, $q>1$.
Clearly, $(C'_n,F'_n)_n$ satisfy (\ref{CF1})--(\ref{CFfin}).
Let $T^{(2)}$ be $(C,F)$-action associated with $(C'_n,F'_n)_n$.

Then by construction $T = T^{(1)} \times T^{(2)}$ is an infinite measure preserving action of $\Gamma$ such that $U_T(g_{k_n}) \to U_S(d)$ as $n\to\infty$.

Now let $S:=\widetilde{T}$ stand for the Poisson suspension of $T$ (Subsection~\ref{PoisSusp_sec}).
Then $S$ is a free probability measure preserving $\Gamma$-action.
Since the mapping $\mc A_\Gamma (X,\mu) \ni T \mapsto \widetilde{T} \in \mc A_\Gamma(\widetilde{X},\widetilde{\mu})$ is continuous, $U_S(g_{k_n}) \to U_S(d)$ as $n \to \infty$.
\end{proof}

\begin{lemma} \label{limIdisc}
For any `good' sequence $(g_n)_{n=1}^\infty$ in $G$ the following subsets are residual in $\mc A_\Gamma$:
\begin{align*}
\mc I_d &:= \{ T \in \mc A_\Gamma \mid U_T(d) \text{ is a limit point of } \{U_T(g_n)\}_{n=1}^\infty \} \text{ for any } d\in G, \text{ and } \\
\mc O &:= \{ T \in \mc A_\Gamma \mid 0 \text{ is a limit point of } \{U_T(g_n)\}_{n=1}^\infty \}.
\end{align*}
\end{lemma}

\begin{proof}
$\mc O$ and $\mc I_d$, $d \in G$, are $G_\delta$ subsets in $\mc A_\Gamma$ by Lemma~\ref{Pol}.
We note that $\mc O$ and $\mc I_d$ are $\text{Aut}(X,\mu)$-invariant.
By \cite[Claim~18]{FW} the $\text{Aut}(X,\mu)$-orbit of any free $\Gamma$-action is dense in $\mc A_\Gamma$.
Therefore, it remains to show that $\mc O$ and $\mc I_d$, $d \in G$, contain at least one free action.
$\mc I_d$, $d \in G$, are non-empty by Lemma~\ref{LemDIlim}.
Each Poisson $\Gamma$-action is free and belongs to $\mc O$ \cite{OW}.
\end{proof}


\begin{proof}[Proof of Theorem~\ref{MainTh} for $G$ is a discrete countable Abelian group]
Let $\Gamma_1 := G \times J_1^{p_1} \rtimes \Z(p_1)$ be the auxiliary group defined above for $G$ and $p_1$.
Fixing a `good' sequence in $G$ and applying Lemmata~\ref{DSdisc} and \ref{limIdisc} we deduce that there is an action $\widetilde{T}_1 \in\mc A_{\Gamma_1}$ such that $T_1 := \widetilde{T}_1 \upharpoonright G$ is weakly mixing, $\mc M(T_1) = \{p_1\}$ and $U_{T_1}(r_{n,1})\to I$, $U_{T_1}(g_{n,1})\to 0$, where $(r_{n,1})_{n=1}^\infty$, $(g_{n,1})_{n=1}^\infty$ are `good' sequences in $G$.

Now let $\Gamma_2 := G \times J_2^{p_2} \rtimes \Z(p_2)$.
By Lemmata~\ref{DSdisc} and \ref{limIdisc} for a generic $\widetilde{T}_2 \in \mc A_{\Gamma_2}$ the restriction $T_2 := \widetilde{T}_2 \upharpoonright G$ satisfies the following conditions:
\begin{enumerate}
\item \label{wm2}
$T_2$ is weakly mixing,
\item \label{hs2}
$\mc M(T_2) = \{p_2\}$,
\item \label{rU2}
$U_{T_2}(d)$ is a limit point of $\{U_{T_2}(r_{n,1})\}_{n=1}^\infty$ for each $d\in G$,
\item \label{r02}
$0$ is a limit point of $\{U_{T_2}(r_{n,1})\}_{n=1}^\infty$,
\item \label{gI2}
$I$ and $0$ are limit points of $\{U_{T_2}(g_{n,1})\}_{n=1}^\infty$.
\end{enumerate}
Thus $T_1 \times T_2$ is weakly mixing with $\mc M(T_1 \times T_2) = \{p_1\}\diamond \{p_2\}$ by Lemma~\ref{LemRSS}.
Moreover, in view of (\ref{rU2}) and (\ref{gI2}), there are subsequences $(r_{n,2})_{n=1}^\infty$ and $(g_{n,2})_{n=1}^\infty$ of $(r_{n,1})_{n=1}^\infty$ and $(g_{n,1})_{n=1}^\infty$ such that $U_{T_1 \times T_2}(r_n)\to I$, $U_{T_1 \times T_2}(g_n)\to 0$ as $n\to\infty$.

Continuing in the same way, we obtain by induction a sequence of weakly mixing $G$-actions $T_i$ such that $\mc M(T_1 \times\cdots\times T_k) = \{p_1\}\diamond\cdots\diamond \{p_k\}$ for any $k>0$.
It follows, that $T:=T_1 \times T_2 \times \cdots$ is weakly mixing with $\mathcal M(T) = \{p_1\}\diamond\{p_2\}\diamond\cdots$.
\end{proof}

\section{Concluding remarks} \label{Concl}

The scheme of the proof also works for the groups of the form $\R^m \times G$
where $G$ is a discrete countable Abelian group and $m>0$. 
For that we need to construct explicitly a `rigid' $\Gamma$-action as in Lemmata~\ref{LemRigid} and \ref{LemDIlim} for $\Gamma = \R^m \times G\times J^p \rtimes Z(p)$.
Indeed, in both of these lemmata the required action was obtained as the product of two $(C,F)$-actions.
Let us say that an element $g\in\Gamma$ is \emph{`good'} if all but the first coordinate of $g$ vanish.
Then the analog of Lemmata~\ref{LemRigid} and \ref{LemDIlim} for sequences of `good' elements can be
easily proved by constructing separately two $(C,F)$-actions: $\R^m$-action as in
Lemma~\ref{LemRigid} and $G\times J^p \rtimes Z(p)$-action as in Lemma~\ref{LemDIlim}.
Moreover, one may mimic the proof of Lemma~\ref{LemDIlim} to extend it for any locally compact second countable Abelian group.
It follows then that the main result is still true for the classes of locally compact second countable Abelian groups considered in \cite{DS}.





Note that our realizations are weakly mixing but not mixing since they are rigid.
The question if there are mixing realizations of considered sets is still open.
In fact, the set of mixing $G$-actions is meager in $\mc A_G$ endowed with the weak topology.
Therefore the weak topology is not suitable to apply the Baire category argument.
In contrast, Tikhonov introduced another (stronger then the weak) topology on $\mc A_{\Z}$ with respect to which the subset of mixing $\Z$-actions is Polish \cite{Ti1}.
Using this topology he proved via `generic' argument the existence of mixing transformations with homogeneous spectrum \cite{Ti2}.
It looks plausible that this approach may be useful for findind mixing realizations of the sets considered in the present paper.

\end{document}